\documentclass[11pt]{amsart}
\usepackage{amsmath,amssymb,amscd,amsfonts}

\newtheorem{thm}{Theorem}[section]
\newtheorem{lem}[thm]{Lemma}

\newtheorem{prop}[thm]{Proposition}

\newtheorem{assu-nota}[thm]{Assumption--Notation}

\theoremstyle{remark}

\newcommand{\into}{\hookrightarrow}

\newcommand{\C}{\mathbb C}
\newcommand{\Z}{\mathbb Z}
\newcommand{\Q}{\mathbb Q}

\newcommand{\pp}{\mathbb P}
\DeclareMathOperator{\Aut}{Aut}
\DeclareMathOperator{\Pic}{Pic}

\DeclareMathOperator{\Hom}{Hom}

\newcommand{\epsi}{\varepsilon}

\newcommand{\ga}{\gamma}

\newcommand{\Ga}{\Gamma}
\newcommand{\De}{\Delta}
\newcommand{\Si}{\Sigma}
\newcommand{\si}{\sigma}

\newcommand{\fie}{\varphi}

\newcommand{\OO}{\mathcal{O}}
\newcommand{\calO}{\mathcal{O}}

\newcommand{\inv}{^{-1}}

\numberwithin{equation}{section}

\begin{document}
\title[Surfaces with $p_g=0$ and $K^2=6$]{The classification of  surfaces with $p_g=0$,
 $K^2=6$ and non birational bicanonical map}

\author{Margarida Mendes Lopes, Rita Pardini}
\date{}
\begin{abstract} Let $S$ be a minimal surface of general type  with $p_g=0$ and $K^2=6$, such
that its bicanonical map $\fie\colon S\to\pp^6$ is not birational. The map $\fie$ is a morphism
of degree $\le 4$ onto a surface. The case of $\deg\fie=4$ is completely characterized in
[Topology, {\bf 40} (5)
 (2001),  977--991] and the present paper completes the
classification of these surfaces. It is proven that the
degree of $\fie$ cannot be equal to 3, and the geometry of
surfaces with
$\deg\fie=2$ is analysed in detail. The last section contains three examples of  such surfaces,
two of which  appear to be new.
\end{abstract}

\maketitle
\section{Introduction}  
 A minimal surface  of general type with
$p_g=0$
 satisfies the inequalities $1\le K^2\le 9$.  It is known that for $K^2\geq 2$ the image of the
bicanonical map $\fie$  is a surface (\cite{deg}) and that for
$K^2\geq 5$  the bicanonical map is always a morphism of degree $\leq 4$ (\cite{marg}).
Furthermore, if $K^2=9$ then $\fie$ is birational, whilst if $K^2=7,8$ then $\deg\fie\le 2$
and those surfaces for which $\deg\fie=2$ can be characterized (see \cite{london1},
\cite{london2}, \cite{doppio}). 

A complete classification of the minimal surfaces  with $p_g=0$, $K^2=6$ and bicanonical map
of degree
$4$ was given in \cite{topology}, namely showing that all such surfaces are  Burniat surfaces
(see
\cite{peters}).

The results in the present paper complete the classification of the minimal  surfaces of
general type with $p_g=0$ and $K^2=6$ for which the bicanonical map  is not birational.  We
prove the following:

 \begin{thm}\label{tdegree}
 Let $S$ be a minimal  surface of general type with $p_g(S)=0$ and $K^2_S=6$ for which the
bicanonical map
$\fie$ is not birational. Then the degree of $\fie$ is either $2$ or $4$ and the image of
$\fie$ is a rational surface.
\end{thm}

 \begin{thm}\label{fibre}
 Let $S$ be a minimal  surface of general type with $p_g(S)=0$ and $K^2_S=6$ for which the
bicanonical map
$\fie$ has degree $2$. Then:
\begin{enumerate}
\item  there is a fibration $f\colon S\to\pp^1$ such that
 the general fibre $F$ of $f$ is hyperelliptic of genus 3 and $f$ has 4 or 5   double fibres;

\item the bicanonical involution of $S$  induces the hyperelliptic involution on
$F$.

\end{enumerate}

\end{thm}

  In Section \ref{examples} it is also shown that both possibilities in (i) do occur. The first
example is due to Inoue (\cite {inoue}) whilst the two others were, to our knowledge, not known
previously. These examples are  constructed as bidouble covers of rational surfaces.

We want to point out  the striking  similarity of the results of Theorem \ref{fibre} to the
case of surfaces with
$p_g=0, K^2=7,8$, for which also   the non-birationality of the bicanonical map implies the
existence of a genus $3$ fibration with multiple fibres  (see \cite{london2}).  
\smallskip 

 Section~\ref{secdegree}  is devoted to proving   Theorem \ref{tdegree}, whilst in
Section~\ref{deg2}
 we prove Theorem \ref{fibre}. Finally, in Section \ref{examples} we describe the examples.

\medskip
\noindent{\bf Acknowledgements.} 
 The present collaboration takes place in the framework of the european contract EAGER, no.
HPRN-CT-2000-00099. It was partially supported by the 1999 Italian P.I.N. ``Geometria sulle
Variet\`a Algebriche'' and by the ``Financiamento Plurianual'' of CMAF. The first author is a
member of CMAF and of the Departamento de Matem\'atica da Faculdade de Ci\^encias da
Universidade de Lisboa and the second author is a member of GNSAGA of CNR. 

\medskip
\noindent{\bf Notation and conventions.} We work over the complex numbers; all varieties are
assumed to be compact and algebraic.

 We do not distinguish between line bundles and divisors on a smooth variety, using the
additive and the multiplicative notation interchangeably.  Linear equivalence of divisors is
denoted by
$\equiv$ and numerical equivalence by $\sim$. The remaining notation is standard in algebraic
geometry.

\section{The degree }\label{secdegree}

In this section we prove  Theorem \ref{tdegree}.
 We start by stating some general facts:
 
\begin{lem}\label{ldegree} Let $S$ be a minimal surface of  general type with $p_g(S)=0$ and
$K^2_S=6$ for which the bicanonical map $\fie$ is not birational. Then the degree of
$\fie$ is less than or equal to $4$ and the image of $\fie$ is a rational surface. 
\end{lem}
\begin{proof} Denote by $\Upsilon$ the image of $\fie$, by $r$ the degree of $\Upsilon$ and by
$d$ the degree of the bicanonical map. Recall that 
$K^2_S=6$ implies that 
$h^0(S,2K_S)=7$ and
$(2K_S)^2=24$. Since $|2K_S|$ is base point free by Reider's theorem (\cite{red})  and
$\Upsilon\subset
\pp^6$ is a non-degenerate surface,  one has the following possibilities for the pair
$(d,r)$:
$(4,6)$, $(3,8)$ and $(2,12)$ and therefore $d\leq 4$. In
the two first cases $\Upsilon$, being a surface
of degree
$\leq 10$ in $\pp^6$ with $p_g=q=0$, is of course a rational surface (cf. 
Lemma 1.4 and Remark 1.5 of
\cite{beau}). For the last case $\Upsilon$ is a rational surface by
Theorem 3 of \cite{xiao2}.

\end{proof}

By Lemma \ref{ldegree}, to prove Theorem \ref{tdegree} we have to exclude the case
$\deg\fie=3$. In this case, as seen in the proof of Lemma \ref{ldegree}, the bicanonical image
is a linearly normal rational  surface $\Upsilon$ of degree $8$ in $\pp^6$. We now  establish
some properties of
$\Upsilon$. Notice that surfaces of degree 8 in $\pp^6$ have been studied classically by
Castelnuovo and later by P. Ionescu (\cite{ionescu}), in the smooth case, and by E. Halanay
(\cite{halanay}), in the normal case. 

\begin{prop}\label{rat1} Let $\Upsilon$ be a linearly normal rational surface of degree
$8$ in $\pp^6$, let $\rho\colon X\to \Upsilon$ be the minimal  desingularization of
$\Upsilon$ and  let $H:=\rho^*\OO_{\Upsilon}(1)$.

Then  $\Upsilon$ has  isolated singularities and one of the following occurs: 
\begin{enumerate}
\item    $-K_X$ is nef and big,  $H=-2K_X$ and $h^0(X,-K_X)=3$;
\item  $X$ has a pencil $|C|$ of rational curves such that $HC=2$;
\item $X$ has a pencil $|C|$ of rational curves such that $HC=3$.
 \end{enumerate}
 \end{prop}
\begin{proof} For the reader's convenience we break the proof into steps:
\smallskip

\noindent{\bf Step 1:} {\em The general hyperplane section of $\Upsilon$ is smooth of genus 3.}

\smallskip Let  $H\in |H|$ be   general and set $H':=\rho(H)$, so that $H\to H'$ is the
normalization map.  The curve
$H'\subset \pp^5$ has degree 8, hence by Castelnuovo's theorem (cf. \cite{Ci}) its 
geometrical  genus
$g(H)$ is lesser  than or equal to
$3$.  Since $q(X)=0$, the restriction map
$H^0(X,H)\to H^0(H,H)$ is surjective  and therefore
$h^0(H,H)=6$.   On the other hand,  Riemann--Roch  gives
$6=h^0(H,H)\geq 8+1-g(H)$, and so
$g(H)\geq 3$. Therefore
$g(H)=3$ and by Castelnuovo's theorem, $H'$ is smooth and $\Upsilon$ has only isolated
singularities.

\medskip
\noindent{\bf Step 2:} {\em The divisor $K_X+H$ is nef.}

\smallskip Since $q(X)=p_g(X)=0$, the restriction map
$H^0(X, K_X+H)\to H^0(H, K_H)$ is an isomorphism. Since $g(H)=3$ by Step 1, it follows  
$h^0(X,K_X+H)=3$.  Write $|K_X+H|=|M|+D$,  where $|M|$ is the moving part and $D$ is the fixed
part. Since 
$|K_X+H|$ cuts out the complete linear system  $|K_H|$, which is free for general  
$H$, necessarily we have $H\Ga=0$ for every component $\Ga$ of $D$. Then by the index theorem
one has
$\De^2<0$ for every divisor $\De$  whose support is contained in the support of $D$.

Assume that $K_X+H$ is not nef and let $\theta$ be an irreducible curve such that
$(K_X+H)\theta<0$. Since 
$\theta M\geq 0$, necessarily
$\theta D<0$ and so $\theta$ is a component of $D$. Therefore $H \theta=0$ and 
$\theta^2<0$. The conditions 
$\theta (K_X+ H)<0$ and
$\theta H=0$ imply $\theta K_X<0$, so that  $\theta$ is a $-1-$curve.  Now, because 
$H\theta=0$, this is a contradiction to the assumption that $\rho\colon X\to \Upsilon$ is the
minimal  desingularization of $\Upsilon$. 

\medskip

 Since
$H^2=8$,  by the adjunction formula we conclude that  
$K_XH=-4$.  We have $MH=MH+DH=(K_X+H)H=4$. By the index theorem this implies $(K_X+H)^2\le 2$
and, by the nefness of $K_X+H$, 
$M^2\le(K_X+H)^2\le 2$. On the other hand $M^2\ge 0$ since $|M|$ has no fixed component. 
Notice that  either $D=0$   or  
$DM>0$, since $K_X+H$ is nef and $D^2<0$ if $D\ne 0$.
\medskip

\noindent{\bf Step 3:} {\em If $M^2=2$, then $D=0$ and $H=-2K_X$. In particular $-K_X$ is nef 
and
$h^0(X,-K_X)=3$.}
 
\smallskip
 In this case we have $2M\sim H$ by the index theorem. Thus $DM=0$ and by the above remark  it
follows that
$D=0$. Furthermore, because $X$ is a rational surface,
$2M\sim H$ implies $2M\equiv H$ and so the equality $K_X+H=M$ yields  $M=-K_X$,  $H=-2K_X$.
\medskip

\noindent{\bf Step 4:} {\em If $|M|$ is composite with a pencil $|C|$, then the general curve
$C$ is rational and 
$HC=2$.}

\smallskip Since $h^0(X, M)=3$, in this case we have $|M|=|2C|$. Thus $2\ge M^2=4C^2$, hence
$M^2=C^2=0$. Since $HM=4$, one has $HC=2$. Since $K_X+H$ is nef by Step 1, we get  $2\ge
(K_X+H)^2\ge 2C(K_X+H)$, namely $K_XC\le -1$. Since the general  $C$ is irreducible and
$C^2=0$, one has $K_XC=-2$ and $C$ is a smooth rational curve.
\medskip 

\noindent{\bf Step 5:} {\em If $M^2=1$, then  there is a pencil $|C|$ on $X$ such that the
general
$C$ is rational and 
$HC=3$.}

\smallskip In this case $\phi_M$ is a birational morphism $X\to \pp^2$ by the proof of Step 4.
Since  a general  curve
$M$ in $|M|$ is smooth rational, we have $-2=K_XM+M^2$ and so
$K_XM=-3$.  From $MH=4$ we conclude $M(K_X+H)=1=M^2$ and so $MD=0$, implying
$D=0$. Now the equalities $1=(K_X+H)^2=K_X^2-8+8$ mean that $K_X^2=1$ and thus $X$ is $\pp^2$
blown-up in 8 points, possibly infinitely near. Let
$E_1,...,E_8$ be the corresponding exceptional divisors. Since $K_X=-3M+E_1+\dots+E_8$, we have
$H=4M-E_1-\dots-E_8$ and there is at least a pencil $|C|$ of rational curves on $X$
(corresponding to the lines in $\pp^2$ passing through one of the blown-up points of
$\pp^2$) such that 
$HC=3$. 
\end{proof}

\medskip

In order to prove Theorem \ref{tdegree} we need also the following technical result.
\begin{lem}\label{h1} Let $S$ be a minimal surface of general type with $p_g(S)=0$, $K_S^2=6$
and let
$|F|$ be a rational pencil on
$S$ such that $F^2=1$, $K_SF=3$. Then $h^1(S, 2K_S-F)=0$.
\end{lem}
\begin{proof} We argue by contradiction. Assume that $h^1(S, 2K_S-F)\neq 0$. Then by the
Riemann-Roch theorem we conclude that $h^0(S, 2K_S-F)\geq 4$. Considering now the long exact
sequence obtained from 
 \begin{equation}\label{seq}
 0\to \calO_S(2K_S-2F)\to\calO_S(2K_S-F)\to\calO_F(2K_S-F)\to 0
 \end{equation}  and using the fact that, by the Riemann-Roch theorem for curves, one has
$h^0(F,\calO_F(2K_S-F))=3$, we see that $h^0(S, 2K_S-2F)\geq 1$. \par

Since $(K_S-F)^2=1$ and  we are assuming that $h^1(S,2K_S-F)\neq 0$, 
 by the Kawamata-Viehweg's vanishing theorem we conclude that $K_S-F$ is not nef. Let 
$\theta$  be  an irreducible curve such that $(K_S-F)\theta<0$. Notice that in particular
$F\theta>0$.  Then for any effective divisor
$G\in |2K_S-2F|$ the curve $\theta$ is a component of $G$ and $G\theta \leq -2$. Write
$G=\theta+A$.  Since $F$ is nef and $FG=4$, we have  $ F\theta\leq 4$  and $FA=4-F\theta\leq 3$.
Furthermore, because
$4=G^2=\theta G+AG$, we have $AG=4-\theta G\geq 6$. 

We now show that this does not occur by examining the various possibilities for $F\theta$. 

If $F\theta=4$, then $K_S\theta\leq 3$ and so, because $\theta $ is an irreducible curve,
$\theta^2\geq -5$. On the other hand
 by the index theorem $A^2<0$, because $F^2=1$ and $FA=0$. Since  $6\leq AG=A^2+ A\theta$, one
has
$A\theta\geq 7$. Then  $-2\geq \theta G=\theta^2+A\theta$ implies
 $\theta^2\leq -9$, a contradiction. 

If $F\theta=3$, then $K\theta\leq 2$, and as above $\theta^2\geq-4$. Since $FA=1$ by the index
theorem
$A^2\leq 1$ and as above we obtain a contradiction to $\theta^2\geq -4$.

If $F\theta=2$, as before we conclude that $\theta^2\geq -3$, $A^2\leq 4$, implying that
$A\theta\geq 2$ which leads us to the same contradiction.

Finally, if $F\theta=1$ then $K_S\theta=0$ and $\theta$ is a $-2-$curve. In this case an easy
calculation shows that $A^2=6$ and that $A\sim K_S$. We can  write the equality of
$\Q-$divisors:
$K_S-F= \frac{1}{2}A+\frac{1}{2}\theta$. Since $\theta$ is a normal crossings divisor and
$\frac{1}{2}A=\frac{1}{2}K_S$ is  nef and big, by the vanishing theorem  of Kawamata-Viehweg we
obtain
$h^1(S, K_S+K_S-F)= 0$, contradicting our assumption.
\end{proof} 

We are finally ready to prove Theorem \ref{tdegree}.

\begin{proof}[Proof of Theorem \ref{tdegree}] In view of Lemma \ref{ldegree} we need to show
that the case 
$\deg \fie=3$ does not occur.
By Lemma \ref{ldegree} and its proof, in this case the bicanonical image $\Upsilon$ is a
rational surface of degree 8 in $\pp^6$.   We prove the theorem  by excluding all the
possibilities for 
$\Upsilon$ described in Proposition \ref{rat1}.

 In the first place notice that case (iii) is trivially impossible. In fact if $\Upsilon$
contains a  pencil  of rational curves of degree $3$, then the pull back $|F|$ of this pencil 
satisfies
$2K_SF=9$, which is impossible.  

 Now we consider   case (ii), i.e. $\Upsilon$ contains a  pencil of rational curves of degree
$2$. This pencil gives rise to a pencil $|F|$ in $S$ such that $2K_SF=6$, i.e. $K_SF=3$. Hence
$F^2\geq 0$ is odd, and so by the index theorem $F^2=1$ and $g(F)=3$. Since
$2K_SF=6$ and  the image of $F$ is a conic, we conclude that the restriction map
$H^0(S, 2K_S)\to H^0(F,\calO_F(2K_S))$ is not surjective, hence  $h^1(S, 2K_S-F)\neq 0$,
contradicting Lemma \ref{h1}.

So we are left with case (i), namely $H=-2K_X$. Consider the Stein factorization
$X\overset{\eta}{\to}
\overline{X}\overset{\nu}{\to}\Upsilon$ of $\rho\colon X\to \Upsilon$. Since $-K_X=\frac{1}{2}H$
is nef,  the map
$\eta \colon X\to \overline{X}$ contracts only $-2-$curves. Hence $\overline{X}$ is a normal
surface whose singularities are rational double points. In particular $\overline{X}$ is
Gorenstein and
$K_X=\eta^*K_{\overline{X}}$. By the normality of $\overline{X}$, the bicanonical map
$\fie\colon S\to\Upsilon$ induces a morphism $\overline{\fie}\colon S\to\overline{X}$ such that
$2K_S=\overline{\fie}^*(-2K_{\overline{X}})$. Hence
$\xi:=\overline{\fie}^*(-K_{\overline{X}})-K_S$ is a non trivial
$2-$torsion element of $\Pic(S)$ and $h^0(S, K_S+\xi)\ge 3$ by Proposition
\ref{rat1} (i). Let $Y\to S$ be the \'etale double cover  given by $\xi$. The standard formulae
for double covers yield:
$$\chi(Y)=2,\quad K^2_Y=12, \quad q(Y)\ge 2.$$
 On the other hand by Corollary
2.2 of
\cite{topology} we have $K^2_Y\ge 16(q(Y)-1)\ge 16$, a contradiction. This completes the proof.
\end{proof}
 
\section{The case of degree 2 }\label{deg2}   This section is devoted to the proof of Theorem
\ref{fibre}.

Throughout all the section we assume that $S$ is a minimal surface of general type with
$p_g(S)=0$ and $K^2_S=6$ such that the bicanonical map $\fie\colon S\to\pp^6$ has degree 2 onto
its image.  We denote by $\si$ the involution of $S$ induced by $\fie$ and by
$\pi\colon S\to\Si:=S/\si$ the quotient map.

\subsection{Preliminaries }\label{prem} 

Given  a smooth projective surface $Y$
 and $k$
 disjoint  nodal curves  $C_1,\dots, C_k$ of $Y$, we define the  binary code $V$ associated to
$C_1,\dots, C_k$. Consider the map 
$\psi\colon
\Z_2^k\to\Pic(Y)/2\Pic(Y)$ defined by \linebreak $(x_1,\dots, x_k)\mapsto\sum x_i[C_i]$, where
$[D]$ denotes the class of a divisor $D$ in\linebreak  $\Pic(Y)/2\Pic(Y)$. We define  $V$  to
be the kernel of
$\psi$  and we denote by $r$ its dimension.  We say that a curve $C_i$ appears in $V$ if there
exists
$v=(x_1,\dots, x_k)\in V$ with $x_i\ne 0$. We denote by $m$ the number of curves $C_i$ appearing
in
$V$. The weight of an element $v=(x_1,\dots, x_k)\in V$ is the number of indices $i$ such that
the coordinate
$x_i$ is non zero. It is easy to show that the weights of the elements of $V$ are divisible by
4.

This situation has been studied in detail in \cite{nodi}. In particular, if we let $G$  be the
abelian group
$\Hom(V,\C^*)$ then there exists a $G-$cover $p\colon Z\to\Si$ branched precisely over the
nodes of
$\Si$ corresponding to the curves that appear in $V$. The numerical invariants of $Z$ can be
computed explicitly in terms of $r$, $m$ and of the numerical invariants of $Y$. It is sometimes
possible to determine $V$ by studying the properties of $Z$, and viceversa. For instance, if
$Y$ is a
 rational surface  with $b_2(Y)\ge 5$ and the number $k$  of disjoint $-2-$curves is the maximum
possible ($=b_2(Y)-2$), then this technique is used in \cite{nodi} to show that the code $V$
is  the code of ``doubly even'' vectors $DE(s)$, where $k=2s$.  We recall  the definition of
$DE(s)$. Given the  code of even vectors $W=\{\sum x_i=0\}\subset \Z_2^{s}$,  
$DE(s)$ is the image of $W$ via the injection $\Z_2^{s}\to\Z_2^{2s}$ defined by
$(x_1,\dots, x_s)\mapsto (x_1, x_1,\dots, x_s, x_s)$.

We are going to study $V$ in the case in which $Y$ is the minimal resolution of the quotient
surface
$\Si$ of $S$ by the bicanonical involution and the $C_i$ are the exceptional curves of
$Y\to\Si$. We will need the following auxiliary result:
 
\begin{lem}\label{lift} Let $G$ be a finite abelian group, let $Y$ be a smooth projective
variety  and let
$p\colon X\to Y$ be a flat
$G-$cover  with building data $L_{\chi}, D_{(H,\psi)}$ (cf.
\cite{ritaabel}). Let
$\ga\in Aut(Y)$ and let $\tau\in Aut(G)$ be such that for every pair $(H,\psi)$ one has
$\ga^*D_{(H,\psi)}=D_{(\tau (H), \psi\circ\tau\inv)}$ and for every $\chi\in G^*$
one has
$\ga^*L_{\chi}\cong L_{\chi\circ\tau\inv}$. 

Then there exists an automorphism $\ga'\colon X\to X$ that lifts $\ga$.
\end{lem}
\begin{proof} Consider the $G-$cover $p'\colon X'\to Y$ obtained from $p\colon X\to Y$ by
taking base change with $\ga\colon Y\to Y$:
\[
\begin{CD} X' @>{\ga_1}>>X\\ @V{p'}VV @VV{p}V\\ Y@>{\ga}>> Y.
\end{CD}
\]

 If we define a new action of
$G$ on $X'$ by letting $g\in G$ act as $\tau(g)$ and we denote by
$p''\colon X''\to Y$ the $G-$cover thus obtained, then by the assumptions $p''$ and $p$ have
the same building data. By the main result of \cite{ritaabel} there exists an isomorphism of
$G-$covers
$\Psi\colon X\to X''$. The map $\Psi$ can also be regarded as an isomorphism
$\Psi\colon X\to X'$ preserving the covering maps to $Y$ (but not the
$G-$action), and the automorphism
$\ga'\colon X\to X$ can be taken to be   the composition $\ga_1\circ \Psi$.
\end{proof}
\bigskip

We return to our surface $S$ and consider  the involution $\si$ of $S$ induced by the
bicanonical map
$\fie$ and the quotient map
$\pi\colon S\to\Si:=S/\si$.

 The fixed locus of $\si$ is the union  of a smooth curve $R$ and of $10$ isolated points $P_1,
\dots,  P_{10}$ (cf. \cite[\S 2]{london2}). We set 
$B:=\pi(R)$ and $Q_i:=\pi(P_i)$, $i=1,\dots, 10$. The surface
$\Si$ is normal and $Q_1,\dots, Q_{10}$   are ordinary double points, which are the only
singularities of
$\Si$. 
%In particular, the singularities of $\Si$ are canonical and the adjunction formula gives
%$K_S=\pi^*K_{\Si}+R$.
  Let $\epsi\colon \hat{S}\to S$ be the blow-up of $S$ at $P_1,
\dots,P_{10}$ and let
$E_i$ be the exceptional curve over $P_i$, $i=1, \dots, 10$. It is easy to check that $\si$
induces an involution
$\hat{\si}$ of
$\hat{S}$ whose fixed locus is the union of $R_0:=\epsi\inv R$ and of $E_1, \dots, E_{10}$.
Denote by
$\hat{\pi}\colon \hat{S}\to Y:=\hat{S}/\hat{\si}$ the projection onto the quotient and set
$B_0:=\hat{\pi}(R_0)$,
$C_i:=\hat{\pi}(E_i)$, $i=1, \dots, 10$. The surface $Y$ is smooth and the $C_i$ are disjoint
$-2-$curves. Denote  by $\eta\colon Y\to
\Si$ the morphism induced by
$\epsi$. The map $\eta$ is the minimal resolution of the singularities of $\Si$ and  there is  a
commutative diagram:
\begin{equation}\label{diagram}
\begin{CD}\ \hat{S}@>\epsi>>S\\ @V\hat{\pi} VV  @VV\pi V\\ Y@>\eta >> \Si
\end{CD}
\end{equation} We recall that by Lemma \ref{ldegree} $Y$ and $\Si$ are rational surfaces. The
map
$\hat{\pi}\colon \hat{S}\to Y$ is a flat double cover branched on $B_0$ and on the $C_i$, hence
it is given by a relation
$2L\equiv B_0+C_1+\dots +C_{10}$. The moving part of the bicanonical system
$|2K_{\hat{S}}|$ is equal to $\hat{\pi}^*|2K_Y+B_0|=\epsi^*|2K_S|$. In particular, since
$|2K_S|$ is free, $|2K_Y+B_0|$ is also free.
\begin{lem}\label{numeri} One has:
\begin{enumerate}
\item $K^2_Y\ge -4$;
\item $L^2+K_YL=-2$, $K_Y^2+K_YL=0$;
\item if $K^2_Y=-4$, then $B_0^2=-4$, $K_YB_0=8$.
\end{enumerate}
\end{lem}
\begin{proof}  Let $\Ga\subset H^2(Y,\C)$ be the subspace generated by the classes of the
curves $C_1,\dots, C_{10}$.  Since the intersection matrix $(C_iC_j)$ is negative definite,
$\Ga$ has dimension 10 and the orthogonal subspace $\Ga^{\perp}$ has dimension
$h^2(Y)-10=-K^2_Y$. Since there is an injection
$\Ga^{\perp}\into H^2(S,\C)$, we get $4=h^2(S)\ge -K^2_Y$, namely $K^2_Y\ge -4$.

By Riemann--Roch and by the standard formulae for double covers, the condition $L^2+K_YL=-2$ is
equivalent to $\chi(S)=1$. The condition  $K^2_Y+K_YL=0$ expresses the fact that the
bicanonical map
$\fie$ factorizes through the cover $\pi\colon S\to\Si$ (cf. \cite{london2}, Proof of Prop.
2.1). Finally, (iii) follows from (ii) and from the relation $2L\equiv B_0+C_1+\dots+ C_{10}$.
\end{proof}
\medskip

\subsection{The proof of Theorem \ref{fibre}.}

We start by analysing the code $V$ associated to the curves
$C_1,\dots, C_{10}$ on $Y$.
\begin{lem}\label{code} In the above setting:
\begin{enumerate}
\item $r:=\dim V\ge 3$, and if  equality holds then $K^2_Y=-4$;
\item $V$ is the code of doubly even vectors (cf. \S \ref{prem}).
\end{enumerate}
\end{lem}
\begin{proof}  Consider the map $\psi\colon \Z_2^{10}\to\Pic(Y)/2\Pic(Y)$ introduced in \S
\ref{prem}. The intersection form on  $\Pic(Y)$ induces a non degenerate $\Z_2-$valued 
bilinear form on
$\Pic(Y)/2\Pic(Y)$ and the image of $\psi$ is a totally isotropic subspace. Hence: $$2\dim Im
\psi\le
\dim_{\Z_2}\Pic(Y)/2\Pic(Y)=10-K^2_Y\le 14,$$ where the last inequality follows by Lemma
\ref{numeri}, (i). Thus the dimension of $V=\ker\psi$ is at least 3, and if it is equal to 3
then $K^2_Y=-4$. This proves (i).

If the number $m$ of curves appearing in $V$ is $\ge 8$, then statement (ii) follows by
\cite[Theorem 3.2]{nodi}. 

So assume that $m<8$. Since $\dim V\ge 3$ and the weights of $V$ are divisible by 4,
necessarily  $\dim V=3$ and $m=7$ (in fact $V$ is the so-called Hamming code).
 
Set $G:=\Hom(V,\C^*)$ ($\cong \Z_2^3$) and consider the $G-$cover
$p\colon Z\to
\Si$ branched on the nodes of $\Si$ corresponding to the curves appearing in $V$ (cf.
\cite[\S 2]{nodi}). By the formulae in \cite[Proposition 2.3]{nodi}, $Z$ is a rational surface
with
$K^2_Z=-32$ and therefore
$b_2(Z)=42$. The surface $Z$ has $24$ nodes.  Let $Z'\to Z$ be the minimal desingularization
and denote by $C_1',\dots, C_{24}'$ the exceptional curves. Arguing as in the proof of (i), one
sees that the code $V'$ associated to the $C_i'$ has dimension $\ge 3$. Since $G$ acts freely
on the set of nodes of $Z$,  the number of curves appearing in $V'$ is divisible by $8$. 
Consider the cover
$q\colon W\to Z$ associated to
$V'$. Again by the formulae of \cite{nodi}, $W$ is an irregular ruled surface. Consider now the
composite map $p':=p\circ q\colon W\to \Si$. 

\smallskip

\noindent {\bf Claim:} {\em $p'$ is a Galois cover.}

\begin{proof}[Proof of claim.]We wish to show that  any given  $\ga\in G$ can be lifted to
$W$.  This  is easier to prove if one works with flat maps, hence we consider the following
diagram, obtained by base change and normalization:
\[
\begin{CD}  W' @>>>W\\ @V{q'}VV @VV{q}V\\ Z'@>>>Z.
\end{CD}
\]  The surface $W'$ is obtained from $W$ by taking the minimal resolution of its singularities
and blowing up
 the ramification points of
$q$. The map
$q'$ is flat.

 We  denote again by $\ga$ the induced automorphism of $Z'$. Clearly,
$\ga$ induces by pull-back an automorphism of the code $V'$ and, dually, an automorphism $\tau$
of the Galois group $G':=Hom(V',\C^*)$ of
$q'$. If we denote by $D_{(H,\psi)}$, $L_{\chi}$ the building data of $p'$, then we 
have $\ga^*D_{(H,\psi)}=D_{(\tau(H),\psi\circ \tau\inv )}$ (cf. Proof of \cite[Prop.
2.1]{nodi}). Since the surface
$Z'$ is simply connected, the line bundles
$L_{\chi}$ associated to the cover $q'$ are determined uniquely by the fundamental relations of
the cover (cf. \cite{ritaabel}) and thus
$\ga^*L_{\chi}\cong L_{\chi\circ\tau \inv}$. So Lemma \ref{lift} applies, $\ga$  can be lifted
to an automorphism $\ga'$ of
$W'$ and it is easy to see that  $\ga'$ induces the required automorphism of $W$, that we
denote again by $\ga'$.
\end{proof}

Now, let $\alpha\colon W\to B$ be the Albanese pencil. The fibration $\alpha$ induces a pencil
$g\colon Z\to \pp^1$ and,  by the Claim, it induces also a pencil of rational curves $h\colon
\Si\to \pp^1$. We have a commutative diagram:
\[
\begin{CD} Z @>p>>\Si\\ @VgVV @VVhV\\ \pp^1@>\bar{p}>>\pp^1.
\end{CD}
\]  The map $\bar{p}$ is a Galois cover with group $G$, since every element of
$G$, having $0-$dimensional fixed locus, acts non trivially on the set of fibres of $g$. On the
other hand this is impossible, since it is well known that $\Aut(\pp^1)$ has no finite subgroup
isomorphic to $\Z_2^3$.  So the case in which  $\dim V=3$ and $m=7$ cannot occur. This
concludes the proof of Lemma \ref{code}.
\end{proof}

\smallskip

\begin{prop}\label{fibration} There exists a  fibration $h\colon \Si\to\pp^1$ such that:
\begin{enumerate} 
\item the general fibre of $h$ is isomorphic to $\pp^1$;
\item $h$ has  at least $r+1$ double fibres.
\end{enumerate} Furthermore, $h$ is uniquely determined  by property (i).
\end{prop}
\begin{proof} The existence of the fibration $h$ follows by \cite[Theorem 3.2]{nodi}, in view
of Lemma \ref{code}. 

To see that $h$ is unique, consider any fibration $h'\colon \Si\to\pp^1$ such that the general
fibre
$H$ of
$h'$ is isomorphic to $\pp^1$.  As explained in the proof of \cite[Theorem 3.2]{nodi}, if 
$p\colon Z\to\Si$ is the cover associated to the code $V$, then $Z$ is an irrational ruled
surface. The inverse image  of $H$ in $Z$ is a disjoint union of smooth rational curves, since
$p$ is unramified in codimension 1. It follows that $h'$ induces on $Z$ a fibration in 
rational curves, which necessarily coincides with the Albanese fibration. Since the cover
$p\colon Z\to \Si$ is canonically associated  to
$\Si$, this shows that $h=h'$.
\end{proof} We denote by $f\colon S\to\pp^1$ the fibration induced by $h\colon\Si\to\pp^1$. The
general fibre
$F$ of $f$ is an hyperelliptic curve and $\si$ induces on $F$ the hyperelliptic involution.
Furthermore the fibration $f$ has at least $r+1$ double fibres corresponding to the $r+1$
double fibres of $h$.
\begin{prop}\label{5fibre} If $r\ge 4$ then $g(F)=3$ and the multiple fibres of $f$ are 5
double fibres.
\end{prop}
\begin{proof}  This proof is very similar to the proof of Theorem 3.2 of  \cite{london2}, but
we include it for the reader's convenience. 

Let
$r\ge 4$. By Proposition
\ref{fibration} and the above remark, the fibration
$f$ has at least $5$ double fibres. 

Suppose that $f$ has at least $6$ double fibres. Let $\psi\colon C\to\pp^1$ be the double cover
branched over the  6 image points of these double fibres of $f$. Note that $C$ is a genus 2
curve. Taking fibre product and  normalization, one gets a commutative diagram:
\begin{equation}\label{diagram2}
\begin{CD}\ X@>\tilde{\psi}>>S\\ @V\tilde{f} VV  @VV f V\\ C@>\psi>> \pp^1
\end{CD}
\end{equation} The map $\tilde{\psi}$ is \'etale, hence $X$ is smooth and we have:
$$\chi(\OO_X)=2\chi(\OO_S)=2,\quad K^2_X=2K^2_S=12.$$  The fibrations $f$ and $\tilde{f}$ have
the same general fibre, which we
 still denote by $F$. Notice   that the genus $g(F)$  of $F$ is odd, since $f$ has double
fibres, and $>1$, since $S$ is of general type. By
\cite{bea},
 we have $$12=K^2_X\ge 8(g(C)-1)(g(F)-1)\ge 16,$$ a contradiction.

Hence $f$ has exactly $5$ double fibres. Then there exists a $\Z_2^2-$cover $\psi\colon  C\to
\pp^1$ branched on the 5 image points of the double fibres of $f$ (cf. \cite{ritaabel},
Proposition 2.1) and again $g(C)=2$. Taking base change and normalization, we obtain a
commutative diagram similar to (\ref{diagram2}). In this case $\tilde{\psi}$ is an \'etale
$\Z_2^2-$cover, hence $K^2_X=24$ and $\chi(\OO_X)=4$. Again \cite{bea} gives
$$24=K^2_X\ge 8(g(C)-1)(g(F)-1)=8(g(F)-1),$$ namely
$g(F)=3$, again because $g(F)$ is odd and $S$ is of general type. 

\end{proof}

\medskip

\begin{proof}[Proof of Theorem \ref{fibre}] By Proposition \ref{5fibre}, to finish the proof 
we have to study the case $r=\dim V= 3$. In this case the argument used in the proof of
Proposition
\ref{5fibre} does not give a bound on the genus of $F$. We obtain such a bound by showing that
the system $|4K_Y+B_0|$ is composed with a pencil $|C|$ of rational curves that satisfy
$CB_0=8$,
$C^2=0$,  so that the pull back of $|C|$ on  $S$ is a genus 3 pencil. Furthermore, we show that
$CC_i=0$ for $i=1,\dots, 10$, and thus   $|C|$ induces a fibration
$h'\colon \Si\to\pp^1$ with rational fibres. By Proposition \ref{fibration}, we have $h=h'$.

\vskip0.2truecm By Lemma \ref{code}, $r=3$ implies 
$K^2_Y=-4$.  We denote by
$H$ a  general element of the linear system $|2K_Y+B_0|$.  Note that, since $|2K_Y+B_0|$ is
free, 
$H$ is a smooth curve. By Lemma
\ref{numeri}, (iii), 
$H^2=12, K_YH=0$ and so $g(H)=7$. We analyse the system $|4K_Y+B_0|=|2K_Y+H|$ by repeated use of
adjunction.  For the reader's convenience we break this analysis into steps:
\smallskip

\noindent{\bf Step 1:} {\em The divisor $K_Y+H$ is nef, $h^0(Y,K_Y+H)=7$ and every curve in
$|K_Y+H|$ is 1-connected.}

\smallskip Since $g(H)=7$, $p_g(Y)=q(Y)=0$, the map
$H^0(Y,K_Y+H)\to H^0(H, K_H)$ is an isomorphism, and so
$h^0(Y, K_Y+H)=7$.  Write $|K_Y+H|=|M|+D$,  where $|M|$ is the moving part and $D$ is the fixed
part. Since 
$|K_Y+H|$ cuts out the complete linear system  $|K_H|$, which is free for general  
$H$, necessarily we have $H\Ga=0$ for every component $\Ga$ of $D$.  Assume that $K_Y+H$ is not
nef and let
$\theta$ be an irreducible curve such that
$(K_Y+H)\theta<0$. Since 
$\theta M\geq 0$, necessarily
$\theta D<0$ and so $\theta$ is a component of $D$. Therefore $H \theta=0$ and  by the  index
theorem
$\theta^2<0$. The conditions 
$\theta (K_Y+ H)<0$ and
$\theta H=0$ imply $\theta K_Y<0$, so that  $\theta$ is a $-1-$curve. We have $\theta B_0=2$,
hence
$\theta$ is not contained in the smooth curve $B_0$  and
 the inverse image $\theta'$  of $\theta$ in $S$ is connected and reduced. Then
$H\theta=0$ implies
$0=\theta'(\hat{\pi}^*H)=\theta'(\epsi^*K_S)$ (see diagram \ref{diagram} for the notation). It
follows that
$p_a(\theta')=0$. Now the Hurwitz formula applied to the double cover $\theta'\to\theta$ gives 
$\theta(B_0+C_1+\dots +C_{10})=2$, i.e. $\theta C_i=0$ for $i=1,
\dots,  10$. Let $Y'$ be the surface obtained from $Y$ by contracting $\theta$ and let
$C_1', \dots, C_{10}'$ be the images of $C_1, \dots, C_{10}$  in $Y'$. The curves $C_1',
\dots, C_{10}'$ are disjoint $-2-$curves and $K^2_{Y'}=-3$. Arguing as in the proof of Lemma
\ref{code} we see that the corresponding code has dimension at least $4$. It is immediate to
check that the code $V$ associated to $C_1,\dots, C_{10}$ has the same dimension, contradicting
the assumption
$r=3$.

So $K_Y+H$ is nef. Furthermore $K_Y+H$ nef  and $(K_Y+H)^2=8>0$ imply that every curve in
$|K_Y+H|$ is 1-connected (see, e.g., \cite[Lemma 2.6]{ml}).
\medskip

\noindent{\bf Step 2:} {\em The general curve $M$ of the linear  system $|K_Y+H|$ is
irreducible, smooth and
$g(M)=3$.}

\smallskip As above write $|K_Y+H|=|M|+D$, where $|M|$ is the moving part and $D$ the fixed
part.  Suppose the general curve in $|M|$ is not irreducible. Then, since $h^0(Y,K_Y+H)=7$ and
$Y$ is regular,  we can write 
$|K_Y+H|=|6G|+D$, where $|G|$ is a pencil with connected fibres such that $GH=2$. Since
$K_Y+H$ is nef and
$(K_Y+H)^2=8$, one has $\frac{8}{6}\geq G(K_Y+H)=6G^2+GD$. Hence we must have $G^2=0$  and 
$GD=1$ (by 1-connectedness). So
$G(K_Y+H)=1$. Since  $GH=2$, we have $GK_Y=-1$, contradicting  
$G^2=0$. 

So the general curve in $|M|$ is  irreducible. Then the image of the map defined by $|M|$ is a
non-degenerate surface in
$\pp^6$, hence of degree $\geq 5$. Since $M^2\leq (K_Y+H)^2=8$,   the linear system $|M|$ can
not have multiple base points, hence the general curve in $|M|$ is smooth.

Now we want to show that $D=0$.  Since $DH=0$ and so $K_YD=MD+D^2$, the number $MD$ is  even. In
addition, we have 
$MD\geq 0$ and  equality holds if and only if $D=0$, by the 1-connectedness of the curves in 
$|K_Y+H|$. 

 Now
$K_Y+H$  nef and
$D H=0$ yield $K_YD\geq 0$,  and so $K_YM\leq -4$, because $K_Y(K_Y+H)=-4$. So $M^2\geq
2g(M)+2$ and so
$h^1(M,M)=0$. We have
$K_YM+M^2=(K_Y+H)M-HM+M^2=2M^2+MD-12$, hence $1-g(M)=-M^2-\frac{1}{2}MD+6$. Since
$h^0(M,M)=6$ and $h^1(M,M)=0$, we obtain by the Riemann-Roch theorem on curves
$6=M^2-M^2-\frac{1}{2}MD+6$, hence
$MD=0$, and so  $D=0$.  So $|K_Y+H|=|M|$, implying also that $g(M)=3$.
\medskip

\noindent{\bf Step 3:} {\em The system $|2K_Y+H|=|K_Y+M|$ is composed with a pencil $|G|$.}

\smallskip Since $M$ is nef and big by Step 1 and Step 2, Kawamata-Viehweg vanishing and
Riemann--Roch give
$h^0(Y, K_Y+M)=3$.

By Lemma \ref{code} (ii), 
$V$ is the code of doubly even vectors. Since $\dim V=3$ by assumption,  we may assume that  the
nodal curves
$C_1$,
$C_2$  do not appear  in the code
$V$ and that $C_3+\dots +C_{10}$ is divisible by 2 in $\Pic(Y)$. Hence the divisors 
$B_0-C_1-C_2$ and
$H-C_1-C_2=2K_Y+B_0-C_1-C_2$ are also divisible by 2 in $\Pic(Y)$ and we may write
$H-C_1-C_2=2\De$ for some $\De\in \Pic(Y)$. Then $\De^2=2$, $K_Y\De=0$ and thus $\chi(\De)=2$
by Riemann--Roch. We have
$(K_Y-\De)H=-6$, hence $h^2(Y,\De)=0$ and $h^0(Y,\De)\ge 2$. Since $Y$ is rational, the
adjunction sequence for $\De$ gives $h^0(Y,K_Y+\De)=h^0(\De,K_{\De})\ge 2$. There is an
inclusion
$|2(K_Y+\De)|\into |K_Y+M|$. Since $\dim |K_Y+M|=2$, the moving part $|G|$ of $|K_Y+\De|$ is a
pencil and $|K_Y+M|$ is composed with $|G|$. 
\medskip

\noindent{\bf Step 4:} {\em The general $G\in |G|$ is isomorphic to $\pp^1$, $G^2=0$ and
$GH=4$.}

\smallskip By Step 3 and its proof,  we can  write $|K_Y+M|=|2G|+D$, where $D$ is the fixed
part. As in the proof of Step 1, we see that $M\Ga=0$ for every component $\Ga$ of $D$, hence
$MG=2$, since $M(K_Y+M)=4$. Since $M^2=8$ and $G$ is nef, the index theorem yields $G^2=0$. Thus
$K_YG=MG-HG=2-HG$ is even by the adjunction formula and it is $\le 0$ since the system $|H|$ is
birational and $G$ moves. If $K_YG=0$ and $HG=2$, then the general $G$ is smooth of genus 1 and
it is mapped by $|H|$ birationally onto a conic, a contradiction. So we have $GH=4$, $GK_Y=-2$
and the general $G$ is smooth rational.

\medskip
\noindent{\bf Step 5:} {\em $|G|$ induces on $\Si$ the fibration $h\colon \Si\to\pp^1$ of
Proposition
\ref{fibration}.}

\smallskip By Step 4 and the uniqueness statement in Proposition \ref{fibration}, it is enough
to show that $|G|$ induces a fibration $\Si\to\pp^1$, namely that $GC_i=0$ for $i=1,\dots, 10$.
Notice that, since
$G(K_Y+M)=G(2G+D)=0$ and $G$ is nef, we have $G\Ga=0$ for every component $\Ga$ of
$D$.  Now assume by contradiction that $GC_i>0$ for some $i$.
 Since $0=C_i(K_Y+M)=2C_iG+C_iD$, this implies that $C_i$ is contained in $D$ and thus $GC_i=0$
by the above remark. 

\medskip
\noindent{\bf Step 6:} {\em The general fibre $F$ of the fibration $f\colon S\to\pp^1$ induced
by
$h\colon
\Si\to\pp^1$ has genus 3.}

\smallskip This follows by applying  the Hurwitz formula to a general $G$, since
$GB_0=G(H-2K_Y)=8$.
\end{proof} 
 
\section{The examples}\label{examples} 

Here we show that Theorem \ref{fibre} is effective, by presenting three examples of  surfaces
satisfying
$K^2=6$, $p_g=0$ and bicanonical map of degree $2$. For one  of the examples the genus
$3$ fibration of Theorem \ref{fibre} has  4 double fibres and  for the other two it has
$5$ double fibres. The first example, due to Inoue (\cite{inoue}, Remark
$6$), is a specialization of a construction of surfaces with $p_g=0$, $K^2=7$ and birational
map of degree 2.  The construction  described here,  as a $\Z_2^2$-cover of a singular rational
surface, is different from the original one, but it has the  advantage of enabling us to
compute the degree of the bicanonical map (cf.
\cite{london1}, Example 1). 
 The other examples are also built as 
$\Z_2^2$-covers of a singular rational surface and are, to our knowledge,  new examples of 
surfaces of general type with these invariants. 
\smallskip 

We start by recalling that, given a smooth surface $Y$,  to define a $\Z_2^2-$cover 
$\pi\colon X\to Y$ one must specify:

 1) divisors $D_1$, $D_2$, $D_3$ of $Y$;

  2) line bundles
$L_1$,
$L_2$ such that $2L_1\equiv D_2+D_3$, $2L_2\equiv D_1+D_3$. 

Notice that if $\Pic(Y)$ has no
$2-$torsion then $L_1$ and $L_2$ are uniquely determined by $D_1$, $D_2$ and $D_3$.  We set
$L_3:=L_1+L_2-D_3$. The divisor $D:=D_1+D_2+D_3$ is the (reduced) branch divisor of $\pi$.
There is a decomposition
$\pi_*\OO_X=\OO_Y\oplus L_1\inv\oplus L_2\inv\oplus L_3\inv$ and the $L_i\inv$ are the
eigenspaces corresponding to the 3 non trivial characters of $\Z_2^2$. We denote by $\chi_i$
the character corresponding to $L_i\inv$ and by $\ga_i$ the generator of $\ker \chi_i$. The
divisor $D_i$ is the image of the divisorial part of the fixed locus of $\ga_i$. The surface
$X$ is smooth if and only if
$D_1$,
$D_2$,
$D_3$ are smooth and
$D$ has normal crossings.
\smallskip

\begin{figure}[ht]
\setlength{\unitlength}{0.5cm}
\centerline{\begin{picture}(24,14)
\thicklines
\put(0,12){\line(1,0){24}}
\put(12,0){\line(0,1){14}}
\put(0,7.5){\line(1,0){24}}
\thinlines
\put(0,13){\line(3,-1){24}}
\put(24,13){\line(-3,-1){24}}
\put(1,14){\line(1,-1){14}}
\put(23,14){\line(-1,-1){14}}
\put(3,12.5){$P_5$}
\put(20,12.5){$P_6$}
\put(7,6.4){$P_1$}
\put(16.3,6.4){$P_3$}
\put(12.5,2.75){$P_2$}
\put(12.2,9.7){$P_4$}
\put(14,12.3){$\De_3$}
\put(12.2,6){$\De_2$}
\put(3.8,7.9){$\De_1$}
\put(0.6,5.75){$S_4$}
\put(23,5.75){$S_3$}
\put(8.7,1){$S_2$}
\put(14.5,1){$S_1$}
\end{picture}}
 \caption{The quadrilateral $P_1P_2P_3P_4$ in $\pp^2$} \label{fig!quad}
 \end{figure}
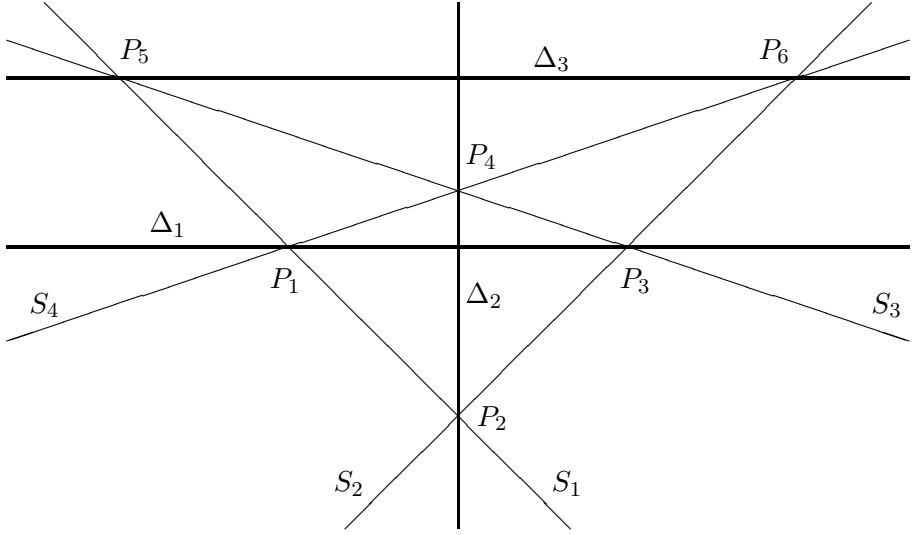

Now consider in $\pp^2$ a quadrilateral $P_1P_2P_3P_4$ (see Figure \ref{fig!quad}). We let
$P_5$ be the intersection point of the lines $P_1P_2$ and $P_3P_4$ and $P_6$ the intersection
point of $P_1P_4$ and
$P_2P_3$. Write $\Si\to\pp^2$ for the blowup of $P_1,\dots,P_6$, and
$e_i$ for the exceptional curves of $\Si$ over $P_i$. Denote by $l$ the pull back of a line.

Write $S_1,\dots,S_4$ for the strict transforms on $\Si$ of the sides
$P_iP_{i+1}$ of the quadrilateral $P_1P_2P_3P_4$ (we take subscripts modulo 4); these are the
only
$-2$-curves of $\Si$. The anticanonical system $|{-}K_{\Si}|$ gives a birational morphism  onto
the symmetric cubic with 4 nodes
$V\subset \pp^3$. This morphism is precisely the contraction of the $-2-$curves of $\Si$ to
canonical singularities. 

If $A\subset\{P_1,\dots,P_6\}$ consists of 4 points no three of which are collinear, then the
linear system of conics through the points of $A$ gives rise to a free pencil on $\Si$; we
denote by $f_1$ the strict transform of a general conic through $P_2P_4P_5P_6$, by $f_2$ that
of a general conic through
$P_1P_3P_5P_6$ and by $f_3$ that of a general conic through
$P_1P_2P_3P_4$.

Finally, we  write $\De_1,\De_2,\De_3$ for the strict transform of the lines $P_1P_3$, $P_2P_4$
and
$P_5P_6$.

The divisors we have introduced satisfy the following relations:
 \begin{enumerate}
  \item $f_i\equiv \De_{i+1}+\De_{i+2}$ for all $i\in \Z_3$;
\item ${-}K_{\Si}\equiv \De_1+\De_2+\De_3\equiv f_1+\De_1\equiv f_2+\De_2\equiv f_3+\De_3$;
 
\item $\De_iS_j=0$ for all $i,j$;
\item $\De_if_j=2\delta_{ij}$ for $1\le i,j\le 3$.
 \end{enumerate}
\medskip
\noindent{\bf Example 1:} Using the above notation, Example 1 of \cite{london1} is obtained by
setting:

1) $D_1=\Delta_1+f_2+ S_1+S_2$,

 $D_2=\Delta_2+f_3$,

$D_3=\Delta_3+f_1+f_1'+S_3+S_4$;

\noindent where $f_1, f_1'\in |f_1|$, $f_2\in |f_2|$, $f_3\in |f_3|$ are general curves, and:

2) $L_1=5l-e_1-2e_2-e_3-3e_4-2e_5-2e_6$,

$L_2=6l-2e_1-2e_2-2e_3-2e_4-3e_5-3e_6$

\noindent and we obtain: $L_3=4l-2e_1-2e_2-2e_3-e_4-e_5-e_6$.

For every $i=1, \dots, 4$, the (set--theoretic) inverse image of $S_i$ in $X$ is the disjoint
union of two $-1$ curves $E_{i1}$, $E_{i2}$:  
 contracting these $8$ exceptional curves on $X$ and contracting the 
$S_i$ on $\Sigma$, one obtains a smooth $\Z_2^2-$cover
$p\colon S\to V$ such that $S$ is of general type with $p_g(S)=0$ and $K^2_S=7$. The
bicanonical map of
$S$ is of degree $2$ and the bicanonical involution coincides with $\ga_1$. The linear system
$|f_1|$ induces  a free pencil $F$ of hyperelliptic curves of genus 3. The bicanonical
involution restricts to the hyperelliptic involution on the general $F$.  The pencil
$|F|$ has 5 double fibres, corresponding to the pull backs of $f_1$, $f_1'$, $\De_2+\De_3$ and
of the two  fibres of $|f_1|$ containing the $-2-$curves.

Now assume that $f_1$, $f_2$ and $f_3$ all  pass through a  general point $P$ and that $f_i$ and
$f_j$ intersect transversely at $P$ for $i\ne j$. In other words, in the terminology of
\cite{quatro} we let the branch locus
$D$ acquire a $(1,1,1)$ point. Denote by $p_0\colon S_0\to V$ the corresponding $\Z_2^2-$cover.
The surface $S_0$ has a singularity of type $\frac{1}{4}(1,1)$ over the image $P'$ of $P$ in
$V$. This singularity can be solved by taking base change with the blow up $\hat{V}\to V$ at  
$P'$ and normalizing. Let $p\colon S\to \hat{V}$ be the cover thus obtained. The exceptional
divisor of $S\to S_0$ is a smooth rational curve with self-intersection $-4$. The surface
$S$ is smooth of general type with $p_g(S)=0$ and $K^2_S=6$ (see \cite{quatro}). A computation
very similar to the one in Example 1 of
\cite{london1} shows that the bicanonical map of $S$ has degree 2 and that the bicanonical
involution coincides with $\ga_1$. As before, the linear system $|f_1|$ induces on $S$ a free
pencil
$|F|$ of hyperelliptic curves of genus 3 such that the bicanonical involution restricts to the
hyperelliptic involution on the general $F$. Now the pencil $|F|$ has 4 double fibres,
corresponding to the pull backs of  $f_1'$, of $\De_2+\De_3$ and of the two  fibres of $|f_1|$
containing the
$-2-$curves. Notice that in this case the pull back of $f_1$ contains with multiplicity 1 the 
exceptional curve of the resolution $S\to S_0$, hence it is not a multiple fibre.
\bigskip

\noindent{\bf Example 2:} With the above notation, consider the point $P_7=\De_2\cap\De_3$ and
denote by $\Si'$ the blow-up of $\Si$ at $P_7$ and by $e_7$ the corresponding exceptional
divisor.

We denote by the same letter the line bundles/divisors on $\Si$ and their pull backs to $\Si'$.
Denote by $\overline\De_2$ and $\overline\De_3$ the strict transforms of $\De_2$ and $\De_3$ and
set:

 1) $D_1=C+ S_1+S_2$,
  
$D_2=f_3$,
 
$D_3=f_1+f_1'+\overline\De_2+\overline\De_3+S_3+S_4$;

\noindent where $f_1, f_1'\in |f_1|$ , $f_3\in |f_3|$ are general curves and $C\in
|f_2+f_3-2e_7|$ is also general;

2) $L_1=5l-e_1-2e_2-e_3-3e_4-2e_5-2e_6-e_7$, and

 $L_2=7l-2e_1-3e_2-2e_3-3e_4-3e_5-3e_6-2e_7$

\noindent and we obtain $L_3=4l-2e_1-2e_2-2e_3-e_4-e_5-e_6-e_7$.

 Since 
$f_2+f_3=4l-2e_1-e_2-2e_3-e_4-e_5-e_6$, it is not difficult to show, for instance by applying a
Cremona transformation centered at $P_1$,
$P_3$, $P_7$, that the general $C\in |f_2+f_3-2e_7|$ is irreducible. Since $p_a(C)=0$, the
general $C$ is also smooth. Thus we obtain a smooth $\Z_2^2-$cover $\pi\colon X\to\Si'$. To
compute the geometric genus of $X$, recall that
$p_g(X)=p_g(\Si')+\sum h^0(\Si', K_{\Si'}+L_i)$ (cf. \cite[Lemma~4.2]{ritaabel}). We have
 \begin{align*}
 K_{\Si'}+L_1 &=2l-e_2-2e_4-e_5-e_6, \\
 K_{\Si'}+L_2 &=4l-e_1-2e_2-e_3-2e_4-2e_5-2e_6-e_7, \\
 K_{\Si'}+L_3 &=l-e_1-e_2-e_3.
 \end{align*}

Clearly both $h^0(\Si',K_{\Si'}+L_1)$ and  $h^0(\Si',K_{\Si'}+L_3)$ vanish.

Now we show that $h^0(\Si',K_{\Si'}+L_2)=0$. Assume otherwise and let $\Ga'\in |K_{\Si'}+L_2|$.
The image
$\Ga$ of $\Ga'$ in $\pp^2$ is a quartic containing $P_1,\dots,P_6, P_7$ which has double
points   at
$P_2,P_4, P_5, P_6$. By Bezout's theorem, the lines in $\pp^2$ corresponding to $S_1$ and $S_2$
are contained in
$\Ga$ and thus
$\Ga'=S_1+S_2+Q'$, where
$Q'$ is the strict transform of a conic $Q$ containing $P_4,P_5,P_6,P_7$ and having a double
point at
$P_4$. But obviously there is no such $Q$ because $P_5,P_6,P_7$ all lie on the line
$\De_3$, which does not contain $P_4$. Hence
$p_g(X)=0$.

For every $i=1,\dots, 4$, the (set--theoretic) inverse image of $S_i$ in $X$ is the disjoint
union of two $-1-$curves $E_{i1}$, $E_{i2}$. Also the inverse image of $\overline{\De}_2$ is the
disjoint union of two $-1-$curves $E_{1}$, $E_{2}$. The bicanonical divisor $2K_X$ is equal to
$\pi^*(2K_{\Sigma'}+D)=\pi^*(-K_{\Sigma'}+f_1+\bar\De_2+S_1+S_2+S_3+S_4)=\pi^*(-K_{\Sigma'}+f_1)+2E_1+2E_2+2\sum
E_{ij}$.

The system $|-K_{\Si'}|$ gives a degree 2 morphism $\Si'\to\pp^2$.  Hence $-K_{\Si'}+f_1$ is
nef and big and it is easy to check that the linear  system 
$|-K_{\Si'}+f_1|$ is birational of (projective) dimension $5$. It follows that the surface $S$
obtained from $X$ by contracting $E_1$, $E_2$ and the $E_{ij}$ is minimal of general type and
the rational map $S\to \Si'$ is composed with the bicanonical map $\fie$ of
$S$. We denote  by the same letter the involutions of $S$ induced by $\ga_1$, $\ga_2$, $\ga_3$.

 Since $2K_X=\pi^*(-K_{\Sigma'}+f_1)+2E_1+2E_2+2\sum E_{ij}$
 one has
$K^2_S=\frac{1}{4}(2K_S)^2=\frac{1}{4}4(-K_{\Sigma'}+f_1)^2=6$.

By the projection formulae for finite flat morphisms, the space $H^0(X, 2K_X)$ decomposes as:
$$
H^0(\Sigma', {-}K_{\Sigma'}+f_1+\overline\De_2+\sum
S_j)\oplus\\
(\oplus_i H^0(\Sigma', {-}K_{\Sigma'}+f_1+\overline\De_2+\sum S_j-L_i)),$$
 where $\Z_2^2$ acts on
$H^0(\Sigma', {-}K_{\Sigma'}+f_1+\overline\De_2+\sum S_j-L_i)$ via the character $\chi_i$. 
Since
$P_2(S)=7$ and $h^0(\Sigma', {-}K_{\Sigma'}+f_1+\overline\De_2+\sum S_j)=h^0(\Sigma',
{-}K_{\Sigma'}+f_1)=6$, it follows that $h^0(\Sigma', {-}K_{\Sigma'}+f_1+\overline\De_2+\sum
S_j-L_i)$ is equal to 1 for one of the indices $i_0\in\{1,2,3\}$ and it is equal to 0 for the
remaining two, so that the bicanonical map has degree 2 and the bicanonical involution is 
$\ga_{i_0}$. A computation shows $h^0(\Sigma', {-}K_{\Sigma'}+f_1+\overline\De_2+\sum
S_j-L_1))=h^0(\Si', e_4+\overline\De_2+S_1+S_2+S_3+S_4)=1$, hence the bicanonical involution of
$S$ coincides with $\ga_1$. The linear system $|f_1|$ induces on $S$ a free pencil
$|F|$ of hyperelliptic curves of genus 3 such that the bicanonical involution restricts to the
hyperelliptic involution on the general $F$. Now the pencil $|F|$ has 5 double fibres,
corresponding to the pull backs of $f_1$, $f_1'$,  of the two  fibres of $|f_1|$ containing the
curves $S_1,\dots, S_4$ and of the  fibre $\overline\De_2+\overline\De_3+2e_7$. Let $2A$ be this
last fibre of
$|F|$ on $S$.  The support of $A$ is the union of an elliptic curve with self--intersection -2,
corresponding to $e_7$, and of a $-2-$curve, corresponding to $\overline\De_3$ (recall that the
inverse image of $\overline\De_2$ in $X$ has been contracted in $S$). The two components of $A$
meet at two points.
\bigskip

\noindent{\bf Example 3:} This is in fact a specialization of  Example 2, obtained by letting
$D_2$ contain the curve
$\overline\De_2$. In this case the cover $\pi\colon X\to\Si'$ is not normal. The normalization
$X'$ of $X$ is again a
$\Z_2^2-$cover of $\Si'$ with branch divisors:

 $D_1=C+\overline\De_2+ S_1+S_2$,

$D_2=\De_1+e_7$,

$D_3=f_1+f_1'+\overline\De_3+S_3+S_4$ (cf. \cite{quatro}).

 The minimal model $S$ of $X'$ is a surface with the same properties as before, but the strict
transform of $\overline\De_2$ is now a $-2-$curve. Furthermore, if we denote again by $2A$ the
reducible double fibre of the  pencil $|F|$ on $S$, then
$A=\theta_1+\theta_2+ 2E$, where $\theta_1$ and $\theta_2$ are disjoint $-2-$curves,
corresponding to
$\overline\De_2$ and $\overline\De_3$ and $E$ is an elliptic curve with $E^2=-1$, corresponding
to
$e_7$. One has
$\theta_1 E=\theta_2 E=1$.
\bigskip

\noindent{\bf Remark:} Notice that for Examples 2 and 3 the divisor  $K_S$ is not ample, in
contrast with the case of Burniat surfaces (cf. \cite{topology}) and of surfaces with
$\deg\fie=2$ and $K^2_S\ge 7$ (cf.
\cite{london2}).

\bigskip

\bigskip

\begin{tabbing} 1749-016 Lisboa, PORTUGALxxxxxxxxx\= 56127 Pisa, ITALY \kill Margarida Mendes
Lopes               \> Rita Pardini\\ CMAF \> Dipartimento di Matematica\\
 Universidade de Lisboa \> Universit\a`a di Pisa \\ Av. Prof. Gama Pinto,  2 \> Via Buonarroti
2\\ 1649-003 Lisboa, PORTUGAL \> 56127 Pisa, ITALY\\ mmlopes@ptmat.ptmat.fc.ul.pt \>
pardini@dm.unipi.it
\end{tabbing}

\end{document}